\documentclass[a4paper,11pt]{article}
\usepackage[english]{babel}
\usepackage[latin1]{inputenc}
\usepackage{amssymb}
\usepackage{amsmath}
\usepackage{amscd}
\usepackage{latexsym}
\usepackage{amsthm}
\usepackage{eucal}
\usepackage{eufrak}
\theoremstyle{plain}
\newtheorem{thm}{Theorem}

\newtheorem{prop}[thm]{Proposition}

\theoremstyle{definition}

\theoremstyle{remark}
\newtheorem{oss}{Remark}
\newtheorem{esem}{Example}

\begin{document}

\linespread{1.1}

\title{Ramification of Quaternion
Algebras over Stable Elliptic Surfaces}

\author{Arvid Perego}

\maketitle

\begin{abstract}
The aim of this work is to study the ramification of quaternion
algebras over the function field of a stable elliptic surface, in
particular over the field of complex numbers. Over number fields
there are criteria for the ramification of quaternion algebras
such as the tame symbol formula. We study how this formula can be
interpreted in a geometrical way, and how the ramification relates
to the geometry of the surface. In particular, we consider stable
complex elliptic surfaces that have four $2-$torsion sections.
\end{abstract}

\section[Introduction]{Introduction}

Let $X$ be an (absolutely) irreducible variety over a field $k$.
Let $BrX$ be the Brauer group of the variety $X$, that is, the
group of equivalence classes of Azumaya algebras on $X$. We have
various descriptions of this group. One of these is given by means
of the étale cohomology of $X$: by \cite{Cald}, Th. 1.1.8, we have
an injective homomorphism
$$BrX\longrightarrow
H^{2}(X,\mathcal{O}_{X}^{*})_{tors}=H^{2}_{\acute{e}t}(X,\mathbb{G}_{m}).$$
We will denote $K:=k(X)$ the function field of $X$, and $BrK$ its
Brauer group, that is, the group of equivalence classes of central
simple $K-$algebras. As shown in \cite{Serr}, Ch 4, §4 and §5, we
have the following description of this group by means of Galois
cohomology:
$$BrK\simeq \lim_{\longrightarrow}H^{2}(Gal(L/K),L^{*}),$$where the
direct limit is taken over the family of finite Galois extensions
$L$ of $K$. $BrK$ is a torsion group. For later use, we will
denote $Br_{n}K$ the $n-$torsion part of $BrK$. We will be
concerned with $Br_{2}K$, and, in particular, with equivalence
classes of quaternion algebras over $K$.

$BrX$ and $BrK$ are deeply related: if $\mathcal{A}$ is an Azumaya
algebra, then, by definition, its stalk $\mathcal{A}_{\eta}$ at
the generic point $\eta$ of $X$ is a central simple $K-$algebra.
This gives a morphism
$$BrX\longrightarrow
BrK,\,\,\,\,\,\,\,\,\,\,[\mathcal{A}]\longrightarrow[\mathcal{A}_{\eta}]$$which
is injective when $X$ is regular (\cite{Arti}). Thus, we can think
of an Azumaya algebra on $X$ as a central simple $K-$algebra which
extends to an Azumaya algebra on $X$. We will make this statement
precise in section 3.

The problem is the following. Let us suppose that we are able to
give an explicit element in $BrK$. Does this element extend to an
Azumaya algebra over $X$? If we have a positive answer to this
question, is the Azumaya algebra obtained in this way non trivial?

Our strategy will be the following. Let
$\pi:X\longrightarrow\mathbb{P}^{1}_{\mathbb{C}}$ be a complex
elliptic surface. By \cite{Sil2}, Ch. III, Rem. 3.1 and Prop. 3.8,
the generic fiber $E_{X}$ of $\pi$ is an elliptic curve over the
field $\mathbb{C}(t)$. In the next section, we describe a way to
construct quaternion algebras over $E_{X}$. Then we study when
these extend to Azumaya algebras over $X$.

\section[The Brauer group of an elliptic curve]{The Brauer group
of an elliptic curve}

We recollect some basic facts about the 2-torsion of the Brauer
group of an elliptic curve $E$ over a field $k$. Let us suppose
that the 2-torsion points of $E$ are $k-$rational, and fix an
isomorphism $E(k)[2]\simeq(\mathbb{Z}/2\mathbb{Z})^{2}$. Let
$$Br^{0}E:=\ker(BrE\stackrel{\epsilon^{*}_{0}}\longrightarrow Brk)$$be
the kernel of the restriction of an Azumaya algebra to the origin
of $E$, where $\epsilon_{0}:Spec\,k\longrightarrow E$ is the point
$0\in E$. We have:

\begin{prop}
\label{prop:iso1} There is a canonical isomorphism $Br^{0}E\simeq
H^{1}(k,E).$
\end{prop}

\proof See \cite{Witt}, Lemma 2.1.\endproof

The Kummer sequence for the 2-torsion part of the elliptic curve
$E$ is $$0\longrightarrow E[2]\longrightarrow E\stackrel{\cdot
2}\longrightarrow E\longrightarrow 0.$$Using the previous
proposition we get the following exact sequence:
$$0\longrightarrow E(k)/2E(k)\stackrel{\delta}\longrightarrow (k^{*}/(k^{*})^{2})^{2}
\stackrel{\gamma}\longrightarrow Br_{2}^{0}E\longrightarrow
0$$(see \cite{Sil1}, Ch. VIII, §2).

Before going on, we fix some notation. A quaternion algebra over a
field $k$ is a four dimensional $k-$algebra $A$ with center $k$,
for which we can find elements $i,j,l\in A$ such that $A=k\oplus
i\cdot k\oplus j\cdot k\oplus l\cdot k$, and such that $ij=-ji=l$.
These elements being fixed, we denote $A$ with the standard symbol
$(a,b)_{k}$, where $a=i^{2}$ and $b=j^{2}$ are elements in
$k^{*}$. When no confusion on this symbol is possible, we will
drop $k$, and we will simply write $(a,b)$ for the corresponding
quaternion algebra. The same symbol will be used for its
equivalence class in the Brauer group of $k$.

\begin{oss}
\label{oss:quat} If $a,b\in k^{*}$, the quaternion algebras
$(a,b)_{k}$, $(b,a)_{k}$, $(a,-ab)_{k}$ are all isomorphic.
Moreover, for any $c\in k^{*}$, $(a,bc^{2})_{k}$ is equivalent (in
$Brk$) to $(a,b)_{k}$, so that $(a,b)_{k}$, $(1/a,b)_{k}$,
$(a,1/b)_{k}$ and $(1/a,1/b)_{k}$ are all equivalent in $Brk$.
Finally, the quaternion algebra $(1,a)_{k}$ is trivial in $Brk$.
\end{oss}

Now fix a Weierstrass equation for $E$:
$$y^{2}=x(x-p)(x-q),$$where $p,q\in k$, then
$R=(0,0)$, $P=(p,0)$ are generators of $E(k)[2]$. Let also
$Q=(q,0)$.

\begin{prop}
\label{prop:deltaegamma}Let $M\in E(k)/2E(k)$. Then $\delta(M)\in
(k^{*}/(k^{*})^{2})^{2}$ is given by:
$$\delta(M)=\begin{cases}
(x(M),x(M)-p)\,\,\,\,\,\,\,\mathrm{if}\,\,\,\,M\neq R,P\\
(q/p,-p)\,\,\,\,\,\,\,\,\,\,\,\,\,\,\,\,\,\,\,\,\,\,\,\,\,\,\,\,\,\,\mathrm{if}\,\,\,\,M=R\\
(p,(p-q)/p)\,\,\,\,\,\,\,\,\,\,\,\,\,\,\,\,\,\,\,\,\mathrm{if}\,\,\,\,M=P\\
(1,1)\,\,\,\,\,\,\,\,\,\,\,\,\,\,\,\,\,\,\,\,\,\,\,\,\,\,\,\,\,\,\,\,\,\,
\,\,\,\,\,\,\,\mathrm{if}\,\,\,\,M=O.
\end{cases}$$
Let $f,g\in k^{*}/(k^{*})^{2}$. Then
$\gamma(f,g)=(x,f)_{k}\otimes(x-p,g)_{k}\in Br^{0}_{2}E$.
\end{prop}

\proof For $\delta$ see \cite{Sil1}, Ch. X, Prop. 1.4. For
$\gamma$ see \cite{Witt}, Proposition 2.2. Note that since we have
chosen $R$ and $P$ as generators of $E(k)/2E(k)$, the equations of
$\delta$ and $\gamma$ are different from what is written in
\cite{Witt} (where the basis is given by $P$ and $Q$).\endproof

Let $\pi:X\longrightarrow\mathbb{P}_{\mathbb{C}}^{1}$ be a complex
elliptic surface which has four 2-torsion sections $s_{R}$,
$s_{P}$, $s_{Q}$ and $O$. The generic fiber $E_{X}$ of $\pi$ has
therefore Weierstrass equation:
$$y^{2}=x(x-p)(x-q),$$with $p,q\in\mathbb{C}(t)$, and we assume that
$s_{R}$, $s_{P}$ and $s_{Q}$ correspond to the points $R=(0,0)$,
$P=(p,0)$, $Q=(q,0)$ (and, clearly, $O$ corresponds to the origin
of $E_{X}$). By Proposition \ref{prop:deltaegamma},
$(x-p,f)_{\mathbb{C}(t)}$, $(x-q,g)_{\mathbb{C}(t)}$,
$(x,h)_{\mathbb{C}(t)}\in Br_{2}^{0}E_{X}$ for every
$f,g,h\in\mathbb{C}(t)$. For example, by Remark \ref{oss:quat} we
have that $(x-p,f)_{\mathbb{C}(t)}=\gamma(1,f)$. Since the
function field $K$ of the surface $X$ is generated by the rational
functions $x,y,t$, we see that $(x-p,f)_{\mathbb{C}(t)}$ (as well
as the others) can be viewed as an element in $BrK$. These will be
the quaternion algebras we will study.

\section[The Brauer group of a smooth surface]{The Brauer group of
a smooth surface}

In this section, $X$ will be a smooth surface over an
algebraically closed field $k$ of characteristic zero. In the
introduction, we saw that there is an injective homomorphism
$BrX\stackrel{\rho}\longrightarrow BrK$. We give a description of
the elements in the image of $\rho$. The main result we will use
is the following:

\begin{thm}
\label{thm:artinmumford}There is an exact sequence:
$$0\longrightarrow BrX\stackrel{\rho}\longrightarrow BrK
\stackrel{\alpha}\longrightarrow\bigoplus_{C\subseteq X}
H^{1}(k(C),\mathbb{Q}/\mathbb{Z}),$$where the direct sum is taken
over the set of irreducible curves in $X$. The map $\alpha$ is
defined in the following way: let $\eta$ be the generic point of
an irreducible curve $C$ on $X$, $\mathcal{O}_{X,\eta}\subseteq K$
the local ring of $X$ in $\eta$, and let $D$ be a central division
algebra over $\mathcal{O}_{X,\eta}$. Then $\alpha(D)$ will be the
cyclic extension $L$ of $k(C)$ associated to $D$.
\end{thm}

\proof See \cite{Arti}, Ch. 3, Th. 1. or \cite{Tann}, Lemme
4.1.\endproof

Let $D=(a,b)$ be a quaternion algebra over $K$, so that $a$ and
$b$ are rational functions on $X$. Let us denote $\alpha(D)_{C}$
the component of $\alpha(D)$ in
$H^{1}(k(C),\mathbb{Q}/\mathbb{Z})$. By Theorem
\ref{thm:artinmumford}, there can only be a finite number of
irreducible curves $C_{1},...,C_{n}\subseteq X$ such that
$\alpha(D)_{C_{i}}\neq 0$. The curve $C_{D}$ given by the union of
these $C_{i}$ will be called \textit{ramification curve} of $D$.
By Theorem \ref{thm:artinmumford}, $D$ extends to an Azumaya
algebra over $X$ if and only if its ramification curve $C_{D}$ is
empty. But how to find this ramification curve?

Let $V$ be the closed subset of $X$ given by the curves on which
$a$ or $b$ have zeroes or poles, and let $U$ be its complement.
Over this open subset we can define the following locally free
sheaf $\mathcal{A}_{a,b}$: if $W\subseteq U$ is an open in $U$, we
define
$$\mathcal{A}_{a,b}(W):=\mathcal{O}_{X}(W)\oplus
i\cdot\mathcal{O}_{X}(W)\oplus j\cdot\mathcal{O}_{X}(W)\oplus
l\cdot\mathcal{O}_{X}(W),$$where $i^{2}=a_{|W}$, $j^{2}=b_{|W}$
and $ij=-ji=l$. Moreover, $\mathcal{O}_{X}(W)$ is the center of
this algebra. Since $a_{|W}$ and $b_{|W}$ are regular functions
without zero on $W$, $\mathcal{A}_{a,b}$ is a sheaf of
$\mathcal{O}_{U}-$modules which has the structure of Azumaya
algebra on $U$. Thus, $\mathcal{A}_{a,b}$ is the extension of the
quaternion algebra $(a,b)$ to $U$, and the ramification curve of
$(a,b)$ is contained in $V$. To determine this curve we have to
study how $\mathcal{A}_{a,b}$ behaves on the irreducible
components of $V$. We use the following well-known tame symbol
formula and, for completeness sake, we give a proof:

\begin{prop}
\label{prop:extension}Let $X$ be a smooth surface over an
algebraically closed field $k$ of characteristic 0, and let
$(a,b)\in Br_{2}(k(X))$ be a quaternion algebra over $k(X)$. Let
$U$ and $V$ be as above, and $C$ be an irreducible component of
$V$. Then $\mathcal{A}_{a,b}$ extends to an Azumaya algebra over
$U\cup C$ if and only if the following rational function on $C$
(called tame symbol):
$$c=(-1)^{v(a)v(b)}a^{v(b)}b^{-v(a)},$$is a square in $k(C)^{*}$,
where $v$ is the valuation for the discrete valuation ring
$\mathcal{O}_{X,\eta}$ at the generic point $\eta$ of $C$.
\end{prop}

\proof First, choose a function $t$ such that $C$ is given by
$t=0$ and $v(t)=1$, so that $\mathcal{O}_{X,\eta}/(t)=k(C)$. Then
we can write $a=\alpha t^{v(a)}$, $b=\beta t^{v(b)}$, where
$v(\alpha)=v(\beta)=0$. By Remark \ref{oss:quat}, it suffices to
study only where $a$ or $b$ have zeroes, that is, we can take
$v(a),v(b)\geq 0$. Since $k$ is algebraically closed, $c$ is a
square in $k(C)^{*}$ if and only if $\alpha^{v(b)}\beta^{-v(a)}$
is a square. This is possible if one of the following is
satisfied:
\begin{enumerate}
\item $v(a)$ and $v(b)$ are even, \item $v(a)$ (resp. $v(b)$) is
even and $\alpha$ (resp. $\beta$) is a square in $k(C)$, \item
$v(a)$ and $v(b)$ are odd, and $\alpha\beta$ is a square in
$k(C)$.
\end{enumerate}
The remaining case, namely $v(a)$ and $v(b)$ are odd and
$\alpha\beta$ is not a square in $k(C)$, is studied in
\cite{Arti}, Ch. 4, Prop. 2, where it is shown that $(a,b)$
ramifies over $C$. So, we study the 3 cases listed above. Write
$v(a)=2m+e$, $v(b)=2n+f$, where $e$ and $f$ can be 0 or 1. The
strategy goes as follows: by purity, it suffices study the problem
locally, so we take $p\in C$ and a neighborhood $U'$ of $p$ in
$X$. We construct an Azumaya algebra $\mathcal{D}$ on $U'$ and an
isomorphism between $\mathcal{D}_{|U\cap U'}$ and
$\mathcal{A}_{a,b|U\cap U'}$. $\mathcal{D}$ is the quaternion
algebra $\mathcal{O}_{U'}\oplus i_{1}\cdot\mathcal{O}_{U'}\oplus
j_{1}\cdot\mathcal{O}_{U'}\oplus i_{1}j_{1}\cdot\mathcal{O}_{U'}$,
where $i_{1}$ and $j_{1}$ are sections of $\mathcal{D}$ over $U'$
such that $i_{1}j_{1}=-j_{1}i_{1}$. Similarly, we have
$\mathcal{A}_{a,b}=\mathcal{O}_{U}\oplus
i\cdot\mathcal{O}_{U}\oplus j\cdot\mathcal{O}_{U}\oplus
ij\cdot\mathcal{O}_{U}$ with $i^{2}=a$ and $j^{2}=b$. So let us
study the 3 cases above:
\begin{enumerate}
\item If $e=f=0$, we take $i_{1}=t^{-m}i$, $j_{1}=t^{-n}j$, and
the isomorphism is given by $i_{1}\mapsto i$, $j_{1}\mapsto j$.
\item If $e=0, f=1$ we take $i_{1}=t^{-m}i$, $j_{1}=t^{-n}j$, so
that $i_{1}^{2}=\alpha$ is a square modulo $t$. We write
$i_{1}^{2}=-z^{2}+ht$, $j_{1}^{2}=\beta t$, where $h=dt^{s}$,
$v(d)=0$. Now take $i_{2}=i_{1}-j_{1}i_{1}$ and $j_{2}=j_{1}$. In
this way,
$$i_{2}^{2}=i_{1}^{2}-(j_{1}i_{1})^{2}=(1+j_{1}^{2})i_{1}^{2}=
(1+\beta t)(-z^{2}+ht)=-z^{2}+tw$$where $w=-\beta z^{2}+h-\beta
th$ is unit. Now take $i_{3}=i_{2}$, $j_{3}=t^{-1}(z-i_{2})j_{2}$,
so that $i_{3}^2=-z^{2}+tw$ and $j_{3}^2=\beta w$. Since these two
have valuation 0, we are done. \item If $e=1$, $f=0$, we take
$i_{1}=t^{-n}j$, $j_{1}=t^{-m}i$ and we go on as in the previous
case. \item If $e=1$, $f=1$, we take $i_{1}=t^{-n-m-1}ij$,
$j_{1}=t^{-n}j$ and we go on as in the second case.
\end{enumerate}\endproof

\section[Criterion for ramification]{Criterion for ramification}

In this section we give a criterion for the ramification of a
quaternion algebra of the form $(x-p,f)$, $(x-q,g)$, or $(x,h)$
over an elliptic surface
$\pi:X\longrightarrow\mathbb{P}^{1}_{\mathbb{C}}$ whose
Weierstrass model
$\widetilde{\pi}:\widetilde{X}\longrightarrow\mathbb{P}^{1}_{\mathbb{C}}$
has equation:
$$y^{2}=x(x-p)(x-q),$$where $p,q\in\mathbb{C}[t]$.
For simplicity's sake, we will only consider stable elliptic
surfaces. We start by giving this criterion, which establishes
when $X$ is stable: we study the fibers of $\widetilde{\pi}$ and
follow the Kodaira classification.

\begin{prop}
\label{prop:stable}Using the same notation as above, let
$p(t)=t^{a}p_{1}(t)$ and $q=t^{b}q_{1}(t)$, with
$p_{1}(0),q_{1}(0)\neq 0$. Moreover, let $p(t)-q(t)=t^{n}h(t)$,
with $h(0)\neq 0$. The fiber $\pi^{-1}(0)$ is stable if and only
if $ab=0$. In this case, the fiber is of type
\begin{enumerate}
\item $I_{2a}$ if $a>0$, $b=0$;\item $I_{2b}$ if $a=0$,
$b>0$;\item $I_{2n}$ if $a=b=0$.
\end{enumerate}
\end{prop}

\proof This can easily be done using the Weierstrass form of the
equation of $\widetilde{X}$, that is $y^{2}=x^{3}+Ax+B$, where $A$
and $B$ depend only on $p$ and $q$. Since $\pi^{-1}(0)$ is stable
if and only if $A(0)B(0)\neq 0$ (see \cite{Mira}, Lecture 2, §3
and Lecture 3, §2, §3), writing out the explicit formulas for $A$
and $B$, we see that the fiber is stable if and only if $ab=0$.
The discriminant of this Weierstrass equation is
$\Delta(t)=t^{2a+2b+2n}p_{1}(t)q_{1}(t)h(t)$, and we get the three
possibilities listed above (see again \cite{Mira}, Lecture 2, §3
and Lecture 3, §2, §3).\endproof

We use Proposition \ref{prop:deltaegamma} to study the
ramification of the quaternion algebras of the form $(x,f)$,
$(x-q,g)$ or $(x-p,h)$. If we change coordinates by $x-q\mapsto x$
or $x-p\mapsto x$, we can always get to a quaternion algebra of
the form $(x,f)$.

We have the following proposition:

\begin{prop}
\label{prop:criterio} Let
$\pi:X\longrightarrow\mathbb{P}_{\mathbb{C}}^{1}$ be a stable
elliptic surface, whose Weierstrass model
$\widetilde{\pi}:\widetilde{X}\longrightarrow\mathbb{P}^{1}_{\mathbb{C}}$
is given by the equation $y^{2}=x(x-p)(x-q)$, with
$p,q\in\mathbb{C}[t]$, $p\neq q$ and $m=max\{deg(p),deg(q)\}$
even. Let $f\in\mathbb{C}(t)$. The quaternion algebra $(x,f)$
extends to an Azumaya algebra over $X$ if and only if for every
zero or pole $t_{0}\in\mathbb{P}^{1}_{\mathbb{C}}$ of $f$, we have
$p(t_{0})q(t_{0})\neq 0$ and $p(t_{0})=q(t_{0})$.
\end{prop}

\proof Using Remark \ref{oss:quat}, it is clear that it suffices
to show the proposition in the case of $f(t)=t$, so that $t_{0}$
is 0 or $\infty\in\mathbb{P}_{\mathbb{C}}^{1}$.

It is clear that $(x,t)$ extends to an Azumaya algebra at least
over the open subset $U$ of $\widetilde{X}$ where $x$ and $t$ do
not have any zero or pole: we blow up any singularity in $U$ and
apply Proposition \ref{prop:extension} (here there is no zero nor
pole of $t$ or $x$). Next, we have to study what happens on
$\widetilde{X}\setminus U$. Since the function $x$ is regular,
$\widetilde{X}\setminus U$ is given by $x=0$ and by the fibers of
$\widetilde{\pi}$ on 0 and $\infty\in\mathbb{P}^{1}$. We start by
studying the fiber over 0.

We begin with the case when $p(0)q(0)\neq 0$ and $p(0)\neq q(0)$,
so that the fiber of $\widetilde{\pi}$ over 0 is smooth. Here, we
don't need to blow up (at least near this fiber). We have
$v(t)=1$, $v(x)=0$. By Proposition \ref{prop:extension}, we get
that $(x,t)$ extends to an Azumaya algebra over
$U\cup\widetilde{\pi}^{-1}(0)$ if and only if $x$ is a square in
the function field of $\widetilde{\pi}^{-1}(0)$. Since the
equation of the fiber is $y^{2}=x(x-p(0))(x-q(0))$, with
$p(0),q(0)\neq 0$ and $p(0)\neq q(0)$, we see that this is not the
case, so that $(x,t)$ ramifies.

As second case, we consider when $p(0)q(0)\neq 0$ and $p(0)=q(0)$,
so $\widetilde{\pi}^{-1}(0)$ has equation $y^{2}=x(x-p(0))^{2}$.
The point $P_{0}=(p(0),0,0)$ is singular. Let
$p(t)-q(t)=t^{n}h(t)$ with $h(0)\neq 0$. We change coordinates in
order to get $P_{0}=(0,0,0)$, so the equation becomes:
$$y=x(x+p)(x+p-q),$$and the quaternion algebra is $(x+p,t)$. We
blow up the surface $\widetilde{X}$ in $P_{0}$. We shall write
down explicit equations for the blow-up: we work in the subvariety
of $\mathbb{A}^{3}\times\mathbb{P}^{2}$ with coordinates
$(x,y,t,(x_{1}:y_{1}:t_{1}))$ defined by the equations
$xy_{1}-x_{1}y=0$, $xt_{1}-x_{1}t=0$ and $yt_{1}-y_{1}t=0$. Where
$x_{1}\neq 0$ we have $y=y_{1}x$ and $t=t_{1}x$, so that the
equation becomes
$$y_{1}^{2}=t_{1}^{n}h(t_{1}x)x^{n}+t_{1}^{n}h(t_{1}x)p(t_{1}x)x^{n-1}+x+p(t_{1}x),$$which
describes a smooth variety. The quaternion algebra is
$(x+p,t_{1}x)$, and we study it over the curves given by $t_{1}=0$
and $x=0$. If $t_{1}=0$, then $y_{1}^{2}=x+p(0)$, so that
$v(t_{1})=1$, $v(x+p(0))=0$ and $x+p$ is a square modulo $t$. If
$x=0$, then $v(x)=1$, but $v(x+p(0))=0$ and $p(0)$ is a square. In
both cases $(x+p,t_{1}x)$ extends to Azumaya algebra on this
curve.

Where $y_{1}\neq 0$ we have $x=x_{1}y$, $t=t_{1}y$, quaternion
algebra $(x_{1}y+p,t_{1}y)$ and equation
$1=x_{1}(x_{1}y+p(t_{1}y))(x_{1}+t_{1}^{n}y^{n-1}h(t_{1}y))$.
Calculations similar to those we did before show that the
quaternion algebra extends to Azumaya algebra also over
$t_{1}y=0$.

There remains the component given by $t_{1}\neq 0$, where we have
$x=x_{1}t$ and $y=y_{1}t$. The quaternion algebra becomes
$(x_{1}t+p,t)$ and the equation is
$$y_{1}^{2}=x_{1}(tx_{1}+p)(x_{1}+t^{n-1}h).$$If $n=1$, this surface
is smooth, while if $n>1$, it has a singular point in
$(0,0,0)$. Let $n=1$: if $t=0$, we have $v(t)=1$, $v(x_{1}t+p)=0$.
By Proposition \ref{prop:extension}, $(x_{1}t+p,t)$ extends to an
Azumaya algebra over this curve if and only if $x_{1}t+p$ is a
square modulo $t$, that is if and only if $p(0)$ is a square in
$\mathbb{C}$, which is the case. This shows that if $n=1$, then
$(x,t)$ extends to an Azumaya algebra over
$\widetilde{\pi}^{-1}(0)$. There remains $n>1$: we must blow up in
$(0,0,0)$ till we get a smooth surface and analyze the quaternion
algebra we obtain as we did for the case $n=1$. By induction on
$n$, we study the ramification of the quaternion algebra
$(xt^{n-1}+p,t)$ over a surface of
equation$$y^{2}=x(xt^{n-1}+p)(x+th).$$With the same notation as
before, if $x_{1}\neq 0$ we get the quaternion algebra
$(x^{n}t_{1}^{n-1}+p(t_{1}x),t_{1}x)$ and equation
$$y_{1}^{2}=(t_{1}^{n-1}x^{n}+p(t_{1}x))(1+t_{1}h(t_{1}x)),$$which
describes again a smooth surface. As before, the quaternion
algebra doesn't ramify over the curve given by $t_{1}x=0$. If
$y_{1}\neq 0$, we get the quaternion algebra
$(x_{1}t_{1}^{n-1}y^{n}+p(t_{1}y),t_{1}y)$ and the equation
$$1=x_{1}(x_{1}t_{1}^{n-1}y^{n}+p(t_{1}y))(x_{1}+t_{1}h(t_{1}y)),$$thus
we don't have any ramification. If $t_{1}\neq 0$, we get the
quaternion algebra $(x_{1}t^{n}+p,t)$ and the equation
$$y_{1}^{2}=x_{1}(x_{1}t^{n}+p)(x_{1}+h),$$which describes a
smooth surface, and we have no ramification.

Finally, we show that if $p(0)q(0)=0$, then $(x,t)$ ramifies over
$t=0$. We can suppose $q(0)=0$ (and $p(0)\neq 0$ by Proposition
\ref{prop:stable}), and we write $q(t)=t^{b}q_{1}(t)$, with
$q_{1}(0)\neq 0$. We use the same notation as before for the
blow-up.

Where $x_{1}\neq 0$ we have $y=y_{1}x$, $t=t_{1}x$, quaternion
algebra $(x,t_{1}x)$ and equation
$$y_{1}^{2}=(x-p(t_{1}x))(1-t_{1}^{b}q_{1}(t_{1}x)x^{b-1})$$which
is smooth. Over the curve $t_{1}=0$ we have $v(t_{1})=1$, $v(x)=0$
so that $(x,t_{1}x)$ extends to an Azumaya algebra over this curve
if and only if $x$ is a square modulo $t_{1}$. By the equation
above, we get $x=y_{1}^{2}+p(0)$, so that $x$ is not a square
modulo $t_{1}$. Since this component does not change under further
blow-ups, $(x,t)$ ramifies over $\widetilde{\pi}^{-1}(0)$.

Next, we study how $(x,f)$ behaves on the fiber
$\widetilde{\pi}^{-1}(\infty)$. Since $(x,t)$ is equivalent to
$(x,1/t)$, we can reduce to the previous case, but the equation is
different: using $s=1/t$, we get
$$y^{2}=x(x-p(1/s))(x-q(1/s)).$$Using $m=max\{deg(p),deg(q)\}$, we
have
$$s^{2m}y^{2}=x(s^{m}x-p_{1}(s))(s^{m}x-p_{2}(s))$$where
$p_{1}(s)=s^{m}p(1/s)$, $q_{1}(s)=s^{m}q_{1}(1/s)$. Since $m$ is
even, we can define $Y=s^{3m/2}y$ and $X=s^{m}x$. Finally, we get
the equation
$$Y^{2}=X(X-p_{1})(X-q_{1}),$$and the quaternion algebra
$(s^{m}X,s)$, which is equivalent to $(X,s)$ since $m$ is even. By
the first part of the proof, we can study this quaternion algebra
on the fiber of 0 and conclude as stated in the proposition.

In order to finish the proof, we shall study how $(x,t)$ behaves
over the curve $\{x=0\}$. The only problem here can be given by
singular points of the surface that lie on this curve and on the
fiber of a point $t_{0}\in\mathbb{P}^{1}$ different from 0 and
$\infty$. We can change coordinates in order to get a quaternion
algebra $(x,t+t_{0})$ and the usual equation. We can use the same
calculations as above for the blow-ups, and look to the transform
of $x=0$. In any case, we see that since $t_{0}\neq 0$ is a square
in $\mathbb{C}$, we do not have any ramification.
\endproof

\begin{oss}
\label{oss:campo} As we saw in the proof, we often make use of the
fact that we are over the field of complex numbers, as here any
number is a square. We can modify the statement for any field $k$
of characteristic 0 asking also that for every zero or pole
$t_{0}$ of $f$, $p(t_{0})$ and $q(t_{0})$ are squares in $k$.
\end{oss}

A natural question is whether $(x,f)$ is non trivial in $BrX$. We
have the following:

\begin{prop}
\label{prop:nontrivial} Let $\pi:X\longrightarrow\mathbb{P}^{1}$
be a stable complex elliptic surface whose Weierstrass model has
equation $y^{2}=x(x-p)(x-q),$ where $p,q\in\mathbb{C}[t]$, $p\neq
q$. Assume that $rk(MW(X))=0$, where $MW(X)$ is the Mordell-Weil
group of $X$. If $f\in\mathbb{C}(t)$ and $(x,f)\in BrX$ is
trivial, then we are in one of the following three cases: $f$ is a
square in the function field of $E_{X}$; $f=\lambda q$ and
$q-p=\mu$; $f=-\lambda q$ and $p=-\mu$, where $\lambda,\mu$ are
squares in the function field of $E_{X}$.
\end{prop}

\proof We can argue as follows: if $(x,f)$ is trivial in $BrX$,
then its restriction to the general fiber $E_{X}$ of $\pi$ is
trivial in $BrE_{X}$. By Proposition \ref{prop:deltaegamma},
$(x,f)=\gamma([f],1)$, where $[f]$ is $f$ modulo squares in the
function field $\mathbb{C}(t)$ of $E_{X}$. It is now clear that if
$(x,f)=0$ in $BrX$, then there must be a $\mathbb{C}(t)-$rational
point $M$ of $E_{X}$ such that $\delta(M)=([f],1)$. Since the
Mordell-Weil group of $X$ is torsion, the only possible $M$ can be
$P$, $Q$, $R$ and the origin $O$ of $E_{X}$. By Proposition
\ref{prop:deltaegamma}, we have only this four possibilities:
\begin{enumerate}
\item $f$ is a square in the function field of $E_{X}$; \item
$q=\lambda f$ and $p=\mu(f-1)$; \item $q=-\lambda f$ and $p=-\mu$;
\item $p=\lambda f$ and $q=\mu f(f-1)$
\end{enumerate}
with $\lambda,\mu$ squares in the function field of $E_{X}$. Since
the surface is stable, by Proposition \ref{prop:stable} we see
that the last case is impossible. The other cases agree with the
hypothesis $(x,f)\in BrX$.
\endproof

We give two examples. The first one, due to Olivier Wittenberg, is
an elliptic K3 surface with Picard number 20. This is described in
\cite{Witt}, Section 3, but here we can study it in a more
geometrical way. The second one was suggested to me by Bert van
Geemen, and shows how important the hypothesis on the annulation
of the rank of the Mordell-Weil group of the elliptic surface is.

\begin{esem}
Let $X$ be the elliptic complex surface over $\mathbb{P}^{1}$
whose Weierstrass model has equation
$$y^{2}=x(x-p)(x-q),$$where $p(t)=3(t+1)^{3}(t-3)$ and
$q(t)=3(t-1)^{3}(t+3)$. It is easy to show that the only singular
fibers are of type $I_{2}$ over $0,3,-3$ and of type $I_{6}$ over
$1,-1,\infty$. Thus $X$ is stable, has Euler number $e(X)=24$
(which implies that $X$ is a K3 surface) and Picard number
$\rho(X)=20$. By the Shioda-Tate formula (see \cite{Mira},
Corollary VII.2.4), $rk(MW(X))=0$. By Proposition
\ref{prop:criterio}, since $p(0)=q(0)=-9$, the quaternion algebra
$(x,t)$ does not ramify over 0. Over $\infty$, following the proof
of Proposition \ref{prop:criterio}, we study the ramification over
0 of the quaternion algebra $(x,s)$ for the surface of equation
$$y^{2}=x(x-p_{1})(x-q_{1}),$$where $p_{1}(s)=3(1+s)^{3}(1-3s)$ and
$q_{1}=3(1-s)^{3}(1+3s)$. Since $p_{1}(0)=q_{1}(0)=3$, $(x,t)\in
BrX$. By Proposition \ref{prop:nontrivial}, $(x,t)$ is non
trivial.

By Remark \ref{oss:campo}, we even get \cite{Witt}, Remark 3.7:
the quaternion algebra $(x,t)\in Br(\mathbb{Q}(X))$ is not
ramified over the field $\mathbb{Q}(i,\sqrt{3})$ (we must have
that $-9$ and $3$ are square in the base field, which is the case
for $\mathbb{Q}(i,\sqrt{3})$).

We can even show that $(x-p,t)$ ramifies, but $(x-p,(t-1)(t+3))$
is in $BrX$ and it is non trivial. To do so, we can write $X=x-p$,
so that the Weierstrass model of $X$ becomes $y^{2}=X(X+p)(X+p-q)$
and the quaternion algebras become $(X,t)$ and $(X,(t-1)(t+3))$.
Now we can apply Propositions \ref{prop:criterio} and
\ref{prop:nontrivial}. Using the same kind of arguments, we get
the same conclusions for $(x-q,t)$ and $(x-q,(t+1)(t-3))$.
\end{esem}

\begin{esem}
Let $X$ be an elliptic complex surface over $\mathbb{P}^{1}$ whose
Weierstrass model has equation
$$y^{2}=x(x-p)(x-q),$$where $p(t)=(t-1)^{2}$ and $q(t)=t^{2}+1$. It is easy
to show that $X$ is rational: the only singular fibers are of type
$I_{4}$ over 1, and of type $I_{2}$ over $i$, $-i$, 0 and
$\infty$, so that the Euler number of $X$ is 12 and Picard number
10 (see \cite{Shio}, Lemma 10.2). By Proposition
\ref{prop:criterio}, we get easily that $(x,t)\in BrX$. Since $X$
is rational, $BrX=0$. In this case, we see that the rank of the
Mordell-Weil group of $X$ is 1, by the Shioda-Tate formula, so
that we cannot apply Proposition \ref{prop:nontrivial}.
\end{esem}

\subsection*{Acknowledgment}
I mostly would like to thank Bert van Geemen for everything he did
for me, teaching, conversations, the time he spent on this work,
which is a conclusion and a revised version of my Tesi di Laurea
at the Università degli Studi di Milano. I wrote a preliminary
version of this work during my staying at the Université de Paris
XI - Orsay, during the spring of 2004, for the Program of
Arithmetic Algebraic Geometry, under the precious and helpful
direction of Jean-Louis Colliot-Thélène and Etienne Fouvry, which
I thank for the possibilities they gave me.

\noindent Arvid Perego, Laboratoire de Mathématiques Jean Leray,
Université de Nantes, 2, rue de la Houssinière, BP 92208, F-44322
Nantes Cedex 03, France

\noindent \textit{E-mail address}:
arvid.perego@math.univ-nantes.fr

\end{document}